\documentclass[10pt,compress]{article}
\usepackage{etex,epic,float}
\usepackage{paralist,booktabs}
\usepackage{comment,paralist,dsfont,setspace,subfigure}
\usepackage{amsmath,amsfonts,mathrsfs,amssymb,wasysym}
\usepackage{epsfig,fancybox,setspace,multirow}
\usepackage[all]{xy}
\usepackage{tikz}
\usetikzlibrary{calc}

\newcommand{\etal}{{\it et al. }}
\graphicspath{{figs//}}

\usepackage{color,subfigure}

\newcommand{\st}{\mathrm{st.} }
\newcommand{\mb}[1]{\boldsymbol{#1}}

\usepackage{hyperref}
\begin{document}

\title % (optional, only for long titles)
{Reducing the Candidate Line List for Practical Integration of Switching into Power System Operation}
%\subtitle{Evidence from India}
\author % (optional, for multiple authors)
{Mohammad Majidi-Qadikolai,  \; Ross Baldick  \\
22nd International Symposium on Mathematical Programming,\\ Pittsburgh, July 2015}
%\institute[Universities Here and There] % (optional)
%{
%%  \inst{1}%
%  Department of Electrical and Computer Engineering\\
%  The University of Teas at Austin
%%  \and
%%  \inst{2}%
%%  Institute of Theoretical Philosophy\\
%%  University There
%}
\date{} % (optional)
{}

%\setupinteractionscreen[option=max]
\maketitle
\begin{abstract}
Optimized operation of the transmission network is one solution to supply extra demand by more efficient use of transmission facilities, and line switching is one main tool to achieve this goal. In this paper, we add extra constraints to OPF formulation to limit the maximum number of switching operations in every hour based on network conditions, and add switching cost in the objective function to represent extra maintenance cost as a result of frequent switching. We also propose an algorithm to remove less important lines for switching in different loading conditions, so OPF with transmission switching will be solved faster for real-time operation. It is applied to a case study with several operation hours. 
\end{abstract}

%  \begin{frame}
%    \frametitle{Outline}
%    \begin{itemize}
%    \item Introduction and Motivation
%    \begin{itemize}
%    \item Why switching?
%    \item Sample related work
%    
%    \end{itemize}
%    \item Overview of  technical issues related to switching
%    \begin{itemize}
%    \item  Transient impact of switching
%      \item Protection system misoperation
%      \item Circuit breaker maintenance
%    \end{itemize}
%    \item Modeling
%    \begin{itemize}
%    \item Proposed algorithm
%    \item Mathematical formulation
%    \end{itemize}
%    \item Case study and numerical results
%    \begin{itemize}
%    \item 13-bus system
%    \item Reduced ERCOT system
%    \end{itemize}
%    \item Conclusion and future work
%    
%    \end{itemize}
%  \end{frame}

  \section*{Nomenclature}
  {\noindent\bfseries Sets and Indices:}\\
  $N_b$: Set of buses with index \textit{i, k, n}\\
  $N_g$: Set of all generators with index \textit{g}\\
  $N_l$: Set of all lines (existing and candidate) with index \textit{l, m}\\
%  $N_o$: Set of all existing lines with index \textit{l, m}\\
%  $N_n$: Set of all candidate lines with index \textit{l, m}\\
%  $N_u$: Set of all existing lines and  selected candidate lines with index \textit{l, m}\\
  $L_k$: Set of lines connected to bus \textit{k} \\
  $G_k$: Set of all generators connected to bus $\textit{k}$ with index ${k}$ \\
  $W_k$: Set of wind generators connected to bus  $\textit{k}$ with index ${k}$ \\
%  $\Phi_l^t$: Set of lines with violated post-contingency flows under outage of line $l$ in load block $t$ \\
%  $N_s^t$: Set of system operation states in load block $t$ with index c ($c=1$ represents normal operating condition)\\
  $T$: The length of time window with index  $\textit{t}$\\
  $|\; \; |$: Size of a set\\
  {\bfseries Parameters:}\\
  $q_i$: Per MWh load curtailment penalty at bus $\it i$ \\
%  $\gamma_i$: Per MWh wind curtailment penalty at bus $\it i$ \\
  $Co_g^t$: Per MWh operation cost for generator $g$ in load block $t$\\ 
  $\zeta_l$: Per switching cost for line $\it l$\\
  $d_i^t$: Demand at bus $\it i$ in load block $\it t$ \\
  $B$: Diagonal matrix of line admittance\\
  $Y$: Reduced admittance matrix (column and row related to reference bus are removed)\\
  $\Psi$: Reduced bus-branch incidence matrix (row related to reference bus is removed)\\
  $P_g^{max,t}$: Maximum capacity of generator $\it g$ in load block $t$\\
  $P_g^{min,t}$: Minimum capacity of generator $\it g$ in load block $t$\\
  $f_l^{max}$: Maximum capacity of line $\it l$\\
  $f_l^{min}$: Minimum capacity of line $\it l$\\
  $M_l$: Big $\it M$ is a large positive number for line $\it l$\\
%  $C^t$: Matrix of contingencies that specifies the status of lines under different contingencies (1 for in service and 0 for out of service lines) for load block $t$ with index $\it c$ \\
  $PTDF$: Power transfer distribution factor\\
  $LODF$: Line outage distribution factor\\
%  $\Gamma_{m,l}^t$: Magnitude of violation in flow of line $m$ in load block $t$ when line $l$ is on outage \\
%  $CII_l^t$: Contingency identification index for outage of line $l$ in load block $t$\\
  %$\delta$: $CII_l^t$ relaxation factor\\ 
%  $\alpha$: Line capacity\\
  {\bfseries Decision Variables:}\\
  $\delta_l$: Binary decision variable for switching line $\it l$\\
  $r_{i,c}^{t}$: MW load curtailment at bus $\it i$ under operation state $c$ in load block $\it t$\\
%  $CW_{i}^{t}$: Aggregated MW wind curtailment at bus $\it i$ in load block~$\it t$\\
  $p_{g}^{t}$: Output power of generator $\it g$ in load block $\it t$\\
  $f_{l,c}^{t}$: Power flow in line $\it l$ under operation state $c$  in load block~$\it t$\\
  $\theta_{i,c}^{t}$: Voltage angle at bus $\it i$ under operation state $c$  in load block $\it t$. 
  $\Delta \theta_{l,c}^{t}$ is voltage angle difference across line $\it l$ under operation state $c$  in load block $\it t$, $\Delta \theta_{l,c}^{t}${= }$\theta_{k,c}^{t}${-}$\theta_{n,c}^{t}$ for line $\it l $ from bus $\it k$ to  bus $\it n$.

  \section{Introduction}
  \subsection{Why Switching?}
  Switching in power system is not a new concept and it used from the early formation of power system until now. During the time, the purpose of transmission switching is expanded. In this paper, we categorized the main purposes of transmission switching into three main categories: 
%  \begin{frame}
%    \frametitle{Why Switching?}
%    \framesubtitle{Different purposes}
    \begin{itemize}
    \item Corrective Action:\\
    Applying switching to clear fault and isolate affected equipment is the primary purpose of installing circuit breakers in power systems.  
    \item Preventive Action:\\
    Reliability switching to back the system to normal condition and prevent load shedding and cascading failures are preventive actions that are applied by system operators usually after occurring a contingency in the network.
    \item Economic Purpose:\\
    In the recent decade, using transmission switching to decrease operation costs by managing congestions/network reconfiguration under normal operating conditions is investigated in literature. Therefore, we should expect more switching in power system in the future.

    \end{itemize}
    %More content goes here
%  \end{frame}

  \subsection{Literature review}
  There are extensive literature on transmission switching with different purposes. In this paper, we briefly review some of them.
%  \begin{frame}
%  \frametitle{Sample Related Work}
  \begin{itemize}
	  \item Switching for improving reliability\\
%		  \begin{itemize}
			     Mazi \etal in~\cite{mazi1986} used transmission switching and bus-bar splitting as preventive actions that can help revealing overloads in lines and preventing cascading failures. In~\cite{shao2005},  Shao \etal incorporated line, bus-bar and shunt element switching into an algorithm for mitigating voltage violations and line overloads in system.
			   Other related papers in this area: \cite{glavitsch1985}\nocite{bacher1986}--\cite{shao2006}
%		  \end{itemize}
	 \item Switching for loss reduction and congestion management\\
%		 \begin{itemize}
			  Fliscounakis \etal in~\cite{flis2007} used piece-wise linear approximation technique to linearize network losses representation in their mixed integer programming formulation. They considered line switching and phase shifter tap changing as tools to manage flows in lines and reduce losses. 
			   Authors in~\cite{granelli2006} used transmission switching and network topology reconfiguration as a tool for congestion management and preventing load shedding. They solved the problem using with both MIP formulation and generic algorithm.
			  There are other options for congestion management in the network like suing FACTS devices and phase shifter transformers. Other related papers in this area: \cite{bacher1988}\nocite{omid2014}\nocite{mjidi_FACTS_2009}\nocite{omid2015}--\cite{majidi_FACTS}
%		 \end{itemize}
  \end{itemize}
  
%  \end{frame}
%  
%    \begin{frame}
%    \frametitle{Sample Related Work}
    \begin{itemize}
    \item Switching for topology optimization/Cost-benefit analysis\\
%    \begin{itemize}
     Ruiz \etal in~\cite{ruiz_reduced_2012} proposed shift factor MIP formulation for topology control. Line opening is modeled with flow cancellation technique. This formulation is compact and its size depends on the size of monitoring and switchable lines. They deployed  this technique for transmission planing as well.
    Hedman \etal in~\cite{hedman_optimal_2008} evaluated the impact of transmission switching on nodal prices, load payments, generation revenues, and flowgate prices. There are several other works by Hedman and his group in this area that are cited in the other related papers section for further reading.
     In~\cite{wu_selection_2013}, Wu \etal Developed a heuristic method for transmission switching that integrated different criteria such as limiting violations of line flows, congestion rents, and production costs. They used Locational Marginal Price (LMP) for decision making in their heuristic method. 
     Other related papers in this area: \cite{ruiz_tractable_2012}\nocite{oneil2005}\nocite{hedman2010}\nocite{hedman2011}--\cite{sahraei2014}
%        \end{itemize}
    \end{itemize}
    
 In this paper we would like to highlight the following concerns for transmission switching for economic purpose:
     \begin{itemize}
     \item Can we implement the selected switching plan in real networks?
     \item Can we solve the problem for large scale systems for real-time operation? 
     \end{itemize}   
%    \end{frame}  
  
  The rest of the paper is organized as follows: in section~\ref{sec:tech}, technical issues related to switching is reviewed.  It is followed by the proposed method  and mathematical formulation in section~\ref{sec:model}. Section~\ref{sec:result} has numerical results and the paper is concluded in section~\ref{sec:conclude}.

\section{An Overview of Technical Issues related to Switching}  \label{sec:tech}  
%\subsection{Technical issues related to switching}
%  \begin{frame}
%  \frametitle{Technical issues related to switching}
%  \begin{itemize}
%  \item Transient impact of switching
%  \item Protection system misoperation
%  \item Circuit breaker maintenance
%  \end{itemize}
%  
%  \end{frame}  
  
%    \begin{frame}
%    \frametitle{Technical issues related to switching}
%    \framesubtitle{Transient impact of switching}
Usually in transmission switching (TS) for economic purposes, steady-state of power system is formulated for TS optimization problem. Based on market  intervals, it is a correct assumption for economic purposes. However, from practical (and reliability) perspective, steady-state analysis is not enough for a TS planning to be implemented. Moreover, switching is not a free action and the extra cost as a result of frequent switching should also be considered. In this section, the transient impact of switching, protection system misoperation, and circuit breaker maintenance is reviewed.

\subsection{Transient impact of switching}
\begin{figure} [H]
    \centering
    \subfigure[Line Switching Voltage Transient]{
    \includegraphics[width=0.75\linewidth]{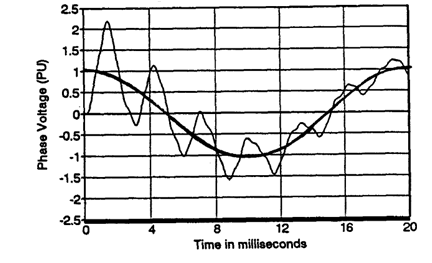}
%    \caption{13 bus network with existing lines, generators and loads.}
  }
   \subfigure[Transformer Energizing Voltage Harmonics
   ]{
      \includegraphics[width=0.75\linewidth]{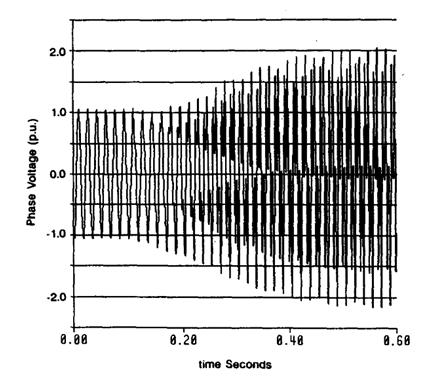}
  %    \caption{13 bus network with existing lines, generators and loads.}
    }
    \caption{Transient impact of switching~\cite{adibi1992}} \label{fig:TS}
    \end{figure}

Opening a transmission line or energizing a transformer will have some transient effects on voltage profile in the system that may trigger  cascading failures. In \figurename{~\ref{fig:TS}} (a), transient over voltage as a result of opening a transmission line is shown. The magnitude of this over-voltage may exceed 2 [P.U.] which is much higher than the accepted voltage deviation (1.05--1.1 [P.U.])during normal operation (steady-state). \figurename{~\ref{fig:TS}} (b) shows the harmonics in voltage as a result of energizing a transformer. These transient phenomena may result in power system protection unnecessary operation that  causes some problems for system reliability. Therefore, they should be considered in switching planning (directly or indirectly).
    
%  \end{frame}  
\subsection{Protection system misoperation}
\begin{figure} [H]
          \centering
          \subfigure[By Category]{
          \includegraphics[width=0.75\linewidth]{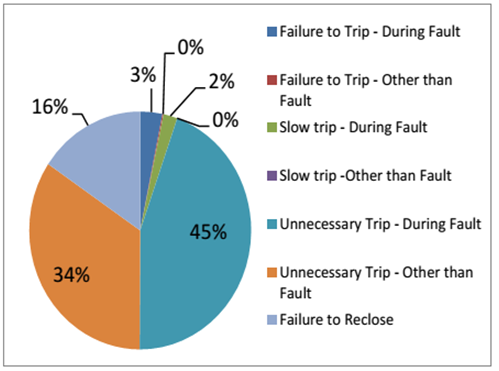}
      %    \caption{13 bus network with existing lines, generators and loads.}
        }
         \subfigure[By Equipment Protected 
         ]{
            \includegraphics[width=0.75\linewidth]{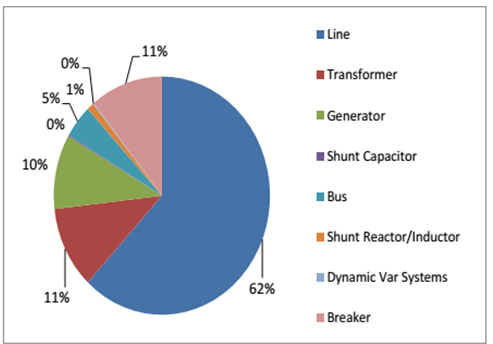}
        %    \caption{13 bus network with existing lines, generators and loads.}
          }
          \caption{Protection System Misoperation Data  ERCOT 2011-2013~\cite{ercot_reliability}} \label{fig:prot}
          \end{figure}

%      \begin{frame}
%      \frametitle{Technical issues related to switching}
%      \framesubtitle{Protection system misoperation}
  Power system protection schemes are mostly designed for static network configurations, and network reconfiguration may result in protection system misoperation. In \figurename{~\ref{fig:prot}}   protection system misoperation in ERCOT from 2011 to 2013 is shown. 
%Although frequent switching for economical purposes are not applied in ERCOT, there are still misoperation in protection system. 
  Adding frequent switching to power system operation will change network configuration more significantly (compared to the case with reliability switching only), therefore it most likely will increase such misoperation if the protection schemes cannot adapt themselves with network reconfiguration (which is the case in many power systems now). 
    
%    \end{frame} 
    \subsection{Circuit breaker maintenance} 
\begin{figure} [H]
        \centering
			 \includegraphics[width=.55\linewidth]{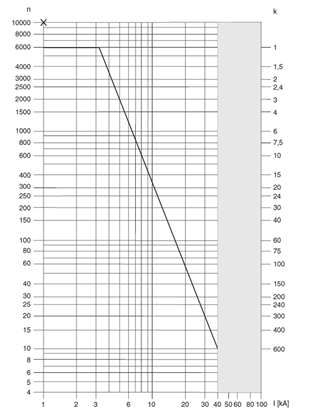}
			 \caption{Number of switching vs. switching current~\cite{cb} \label{fig:cb}
			 }
        \end{figure}

% \begin{frame}
%      \frametitle{Technical issues related to switching}
%      \framesubtitle{Citcuit breaker maintenance}
 As circuit breakers have a limit on the total number of switching before each maintenance, its cost should be considered in switching planning. A SIEMENS circuit breaker switching curve is shown in \figurename{~\ref{fig:cb}}. Based on this figure, total number of switching before  maintenance is 6000 if the breaker opens a circuit under normal condition, and this number will decrease by increasing the current that should be cut by the beaker. For example, if this breaker opens 40KA short circuit current for 10 times it will need  maintenance services. Equation~\eqref{NS} shows how to calculate the number of remaining switching for this breaker depending on the current ($I_x$). $k_i$ and $k_x$ can be found from \figurename{~\ref{fig:cb}} based on current $I_i$ and $I_x$ respectively.
%    \begin{columns}
%        \begin{column}{0.55\textwidth}
%            \begin{itemize}
                
                Number of remaining switching before next  inspection/maintenance \cite{cb}:
                
                \begin{align}
                n_x= \frac{6000-\sum\limits_{i=1}^{m}(n_i \times k_i)}{k_x} \label{NS}
                \end{align}

                where:\\
                $k_i$: weighting factor for $I$\\
                $k_x$: weighting factor for $I_x$\\
                $n_i$: number of performed interruptions at $I$\\
                $n_x$: number of permissible interruptions at $I_k$
              
%            \end{itemize}
%        \end{column}
%        \begin{column}{0.45\textwidth}
          
%        \end{column}
%    \end{columns}

%    \end{frame} 
    In summary, technical issues related to transmission switching should also be considered as a part of planning for switching to be practical.
    
%     \begin{frame}
%          \frametitle{Technical issues related to switching}
%          \framesubtitle{Any impact on swiching planning?}
%          
%          Practical issues on applying switching:
%          \begin{itemize}
%          \item Switching transient and its impact on power system transient stability
%          \item Impact of switching on protection schemes and relay settings
%          \item Direct and indirect costs related to switching
%          \end{itemize}
%          
%          Recent work in this area:\
%          Kory W. Hedman, {Flexible Transmission Decision Support}, PSERC webinar, April, 21, 2015.
%          
%          
%          
%         \end{frame}      
    
  \section{Modeling and Formulation} \label{sec:model}
    
      \subsection{Proposed Algorithm}
%  \begin{frame}
%  \frametitle{Proposed Algorithm}
%  \framesubtitle{Added constraints}
In this paper we have added the following constraints to the classic switching planning problem to make it more practical:
  \begin{itemize}
  \item Adding switching costs to the objective function.
  \item Considering a time window (\textbf{T hours}) to plan for switching rather than a single hour in real-time.
  \item Limiting the number of times that a line can be switched in the planning time window.
  \item Limiting total number of switching  per hour in the network.
  \end{itemize}

%  \end{frame}    
%  
%  \begin{frame}
%    \frametitle{Proposed algorithm}
%    \framesubtitle{Steps for reducing candidate lines lisst for switching}
    As running transient analysis and evaluating the impact of network reconfiguration on protection schemes are computationally expensive, we added last two constraints to integrate expert knowledge into our optimization formulation for limiting the number of switching in a way to be applicable in the system.
    
    Solving switching optimization problem with extra constraints is challenging for large scale power systems especially that it should be run in real-time. We proposed a heuristic method to decrease the search space for   making decision about switching. The proposed algorithm is summarized in the following steps:
    \begin{description}
    \item[Step~1]: Solve OPF for next T hours
    \item[Step~2]: Create Monitored Lines List (MLL)\\
    MLL includes lines that their loading without switching is more than $\alpha\%$.
    \item[Step~3]: Calculate LODF for MLL for all closed lines
    \item[Step~4]: Reduce switching lines list\\
    In this step, lines that their opening will cause overload in lines in MLL will be removed from switching list, as they will have negative impact on monitored lines.
    \item[Step~5] Solve OPF with reduced switching list
       
    \end{description}
    
  As mentioned above, this method is heuristic, as we do not consider multiple line switching in steps 3 and 4. Therefore we cannot guarantee optimality, and the answer will be sub-optimal. However, reducing the switching lines list will reduce computational time. 
    
%  \end{frame} 
  \subsection{Mathematical Formulation}
  Transmission switching optimization problem formulation is given in~\eqref{obj1}--\eqref{tt}. The objective function~\eqref{obj1} includes load shedding penalty cost, electricity generation cost and cost related to switching. In this formulation we penalize breakers operation for both line opening and closing, but it is possible to limit it to line opening if it is preferred. Equations~\eqref{power.balance}--\eqref{gen.cap} represent standard constraints for optimal power flow with line switching binary variable $\delta$. Equations~\eqref{ll} and~\eqref{tt} are two new constraints that are added to integrate expert knowledge into the optimization problem. Equation~\eqref{ll} limits total number of switching for a line during next T hours, and equation~\eqref{tt} limits total number of switching in the system in every hour. The right hand side of these equations will be set by system operators based on the network configuration and loading condition.  
%  \begin{frame}
%  \frametitle{Mathematical Formulation}
%  \framesubtitle{OPF with switching}
% {\footnotesize
 \begin{align}
    Z^{*}{=}&\min\limits_{\mb{p},\mb{\theta},\mb{r},\mb{f}, \mb{\delta}} \sum\limits_T \big[\sum\limits_{N_b}  q_i r_{i}^t {+}\sum\limits_{Ng} Co_g^t p_g^t + \sum\limits_{N_l}\zeta_l |\delta^t_l-\delta^{t-1}_l| \big]\label{obj1} \\
  	\st\quad & {-}\sum\limits_{L_k} f_{l}^{t}{+} \sum\limits_{G_k}p^{t}_{g}{+}r_{k}^{t}{=}d_k^t,\hfill \forall t,k \label{power.balance}\\ 
 	& {-}M_l(1-\delta^t_l)\leq f_{l}^{t}{-}B_{l,l} \Delta\theta_{l}^{t},\hfill \forall t,l \label{lower.flow.1}\\ 
 	& M_l(1-\delta^t_l)\geq f_{l}^{t}{-}B_{l,l} \Delta\theta_{l}^{t},\hfill \forall t,l \label{upper.flow.1}                     \\
 %	\end{align}
 %\begin{align}
 %& CW_{i}^{t}\geq \sum\limits_{W_k} (P_{g}^{max,t}-p_{g}^{t}), \hfill \forall t,i \label{wind.curtail.1}\\
 	   & (\delta^t_l)f_{l}^{min} \leq f_{l}^{t} \leq f_{l}^{max}(\delta^t_l), \hfill \forall t,l \label{line.cap.1}               \\
 %	   \end{align}
 % \begin{align} 
  	  & 0\leq r_{i}^{t} \leq d_i^t, \hfill \forall t,i \label{load.curtail}\\ 
 	  &{-}\frac{\pi}{2} \leq \theta_{i}^{t} \leq \frac {\pi}{2}, \hfill \forall t,i \label{voltage.angle}\\ 
 %	  & CW_{i}^{t}\geq 0, \hfill \forall t,i \label{wind.curtail.2}\\
% 	& -R_g^{Down} \leq p_g^t-p_g^{t-1} \leq R_g^{Up}, \,\forall t,g \label{ramp1}\\
 	  &  P_g^{min,t}\leq p_{g}^{t}\leq P_g^{max,t}, \hfill \forall t,g \label{gen.cap} \\
 	  & \sum\limits_T |\delta^t_l-\delta^{t-1}_l| \le H1_l, \forall l \label{ll}\\
 	  &\sum\limits_{N_l} |\delta^t_l-\delta^{t-1}_l| \le H2^t, \forall t \label{tt}
 %	  & x_l{=}1, \hfill \forall l \in N_o \label{decision.1} \\ 
 %	  &  x_l \in \{0,1\}, \hfill \forall l \in N_l \label{x.value}  
 	\end{align} 
% }
 
%  \end{frame} 

%  \begin{frame}
%  \frametitle{Mathematical formulation}
%  \framesubtitle{MLL, LODF and SLL formulation~\nocite{majidiSDCL}}
%  {\footnotesize
  Equations~\eqref{Viol} and \eqref{MLL} are used to create monitored lines list. MLL contains lines with more than $\alpha\%$ loading. 
  \begin{align}
  & Viol_l^{t}=f_l^t- \alpha\times f_l^{max}, \forall l,t \label{Viol}\\
  & MLL=\{l \in N_l\,|\,Viol_l^t > 0,\forall t \} \label{MLL}
  \end{align}

  To calculate LDOF and post switching line flows, equations~\eqref{PTDF_ll}--\eqref{f_m} are used~\cite{majidiSDCL}.
\begin{align}
& PTDF_{l,l}=B_{l,l} \Psi_l^T [Y]^{-1}  \Psi_l \label{PTDF_ll}\\
& PTDF_{m,l}=B_{m,m} \Psi_m^T [Y]^{-1}  \Psi_l \label{PTDF_ml}\\
& LODF_{m,l}=PTDF_{m,l} (1-PTDF_{l,l})^{-1} \label{LODF}\\
& f_{m,l}^t=f_m^t+LODF_{m,l} f_l^t \label{f_m}
\end{align}

Switching Lines List $SLL$ is reduced by using equations \eqref{CLLE} and \eqref{CLLU}.

\begin{align}
%& PTDF_{l,k}=(\chi_{ik}-\chi_{nk})\,{B_{l,l}} \label{PTDF_ll}\\
%\end{align}
%\begin{align}
%& LCDF_{m,l}=\frac{PTDF_{m,i}^{IN}-PTDF_{m,i}^{OUT}}{PTDF_{l,i}^{OUT}} \label{LCDF}\\
& SLL_E=\{l\in N_l\,|\, f_{m,l}^t-f_m^{max} \geq 0, \forall m \in MLL\} \label{CLLE}\\
& SLL_u=SLL_o\setminus SLL_E \label{CLLU}
%& -1 \leq \beta \leq 0
\end{align}
%}

%\end{frame}
\section{Case study and numerical results} \label{sec:result}
All illustrated results in this section have been obtained from a personal computer with 2.0-GHz CPU using MATLAB R2014a~\cite{MATLAB:2014} and YALMIP R20140221 package~\cite{YALMIP} as a modeling language and GUROBI 5.6~\cite{gurobi} as the solver. Two different case studies consisting of 13-bus system and reduced ERCOT network with 317-bus are considered. 
%    \begin{frame}
%    \frametitle{Case study and numerical results}
%    The proposed method is run for these two case studies:
%    \begin{itemize}
%    \item 13-bus system
%    \item Reduced ERCOT system
%    \end{itemize}
%    \end{frame}
\subsection{13-bus system}
\begin{figure}
        \centering
			 \includegraphics[width=.65\linewidth]{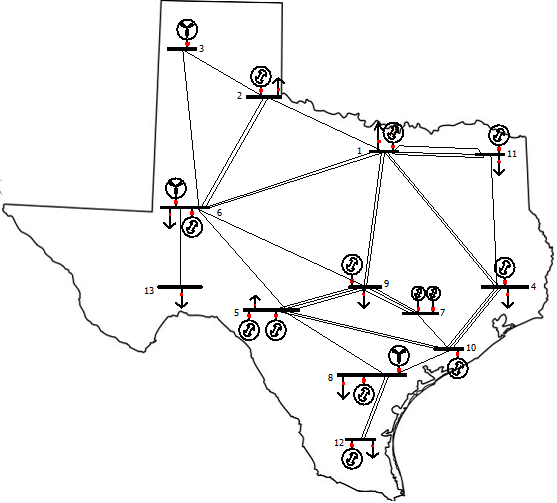}
			 \caption{13-bus system~\cite{majidiSDCL}\label{fig:base.case}
			 }
        \end{figure}

%\begin{frame}
%\frametitle{Case study and numerical results}
%\framesubtitle{13-bus system~\cite{majidiSDCL}}
%
%\begin{columns}
%        \begin{column}{0.55\textwidth}
This 13-bus system is a simplified version of the ERCOT network that is developed for educational purposes (see \figurename~\ref{fig:base.case}). This case study has 13 buses, 33 branches, 16 power plants, and 9 load centers~\cite{majidiSDCL}.
The parameters are set as following:

            \begin{itemize}
               \item Next 5 hours is considered for switching planning ($T=5$).
                \item $H1_l=2$.
                \item $H2^t=4$.
                \item $|N_l|=34$.
                \item $\alpha=50\%$.

            \end{itemize}
%        \end{column}
%        \begin{column}{0.45\textwidth}
%          
%        \end{column}
%    \end{columns}
%\end{frame}

\begin{figure}
        \centering
			 \includegraphics[width=.45\linewidth]{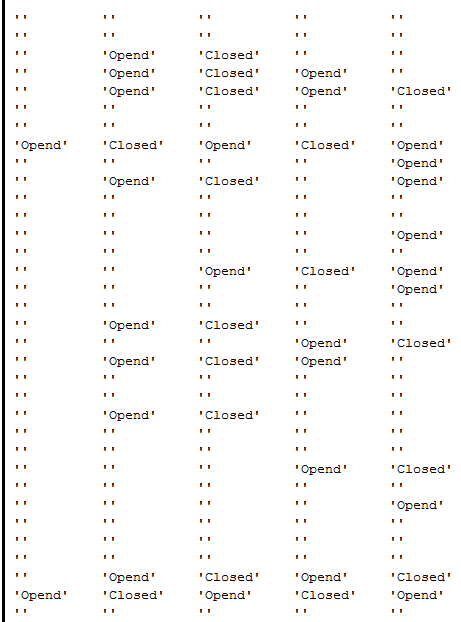}
			 \caption{Switching plan for next 5~hours \label{fig:NC}
			 }
        \end{figure}

%\begin{frame}
%\frametitle{Case study and numerical results}
%\framesubtitle{Switching without any extra constraint}
%\begin{columns}
%        \begin{column}{0.55\textwidth}
In the first step, transmission switching (TS) is solved without any extra constraints defined in this paper (classic TS). The selected TS plan is shown in \figurename{~\ref{fig:NC}}. Operation cost saving as a result of switching is 0.12\% compared to the base case without any switching.  Without any extra constrains, the system operator should switch 44 times during next 5 hours in a network with 34 lines. Lines 8 and 33 are switched 5 times during 5 hours. 12 times switching in hour 5 will significantly change system configuration and may affect protection system performance that makes such switching plans impractical.

%           Saving as a result of switching: 0.12\% \\
%           How about reliability?
%        \end{column}
%        \begin{column}{0.45\textwidth}
          
%        \end{column}
%    \end{columns}
%
%\end{frame}

\begin{figure}
          \centering
          \subfigure[]{
          \includegraphics[width=0.45\linewidth]{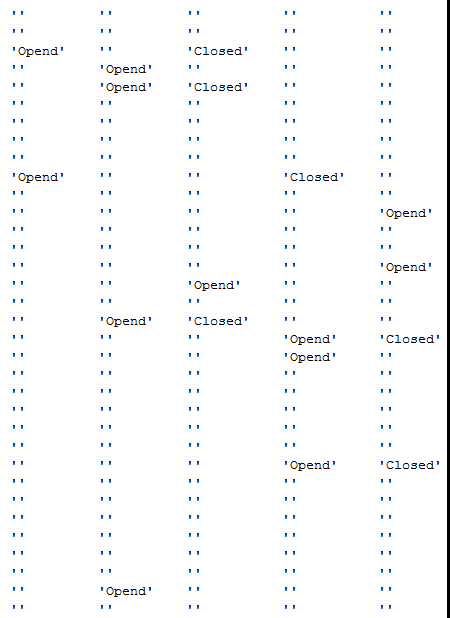}
      %    \caption{13 bus network with existing lines, generators and loads.}
        }
         \subfigure [ 
         ]{
            \includegraphics[width=0.45\linewidth]{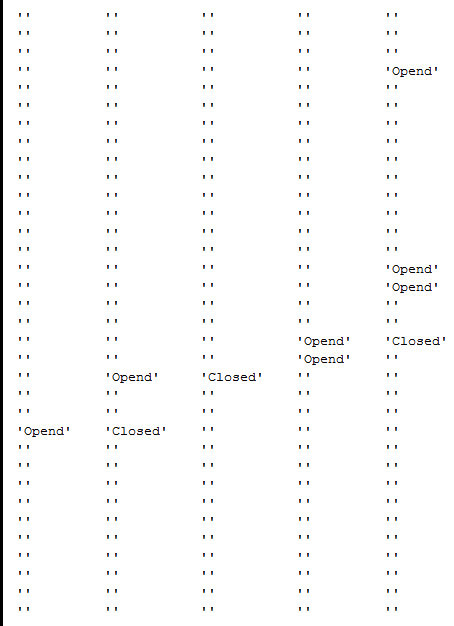}
        %    \caption{13 bus network with existing lines, generators and loads.}
          }
          \caption{Switching plan with (a) and without (b) extra constrains } \label{fig:SP}
          \end{figure}

%\end{frame}
By adding extra constraints (still not using MLL and SLL), the total number of switching will reduce to 18 (shows 144\% switching reduction). However, this less switching will result in less saving and some extra operation costs compared to the TS without any extra constraint. For this case study, this extra cost is  0.021\%. \figurename{~\ref{fig:SP}} (a) shows the switching plan for next 5 hours after adding extra constraints. As we didn't use our proposed heuristic method until now, this TS is optimal.
%\begin{frame}
%\frametitle{Case study and numerical results}
%\framesubtitle{Impact of adding new constraints on switchign plan}
%{\footnotesize
%    Extra operation costs: 0.021\%\\ 
%    Total number of switching: 18\\
%     Switching reduction: 144\%
 
%}
%
%\end{frame}

In the next step, we create MLL and SLL based on equations in section~\ref{sec:model}. As shown in Table~\ref{tab:Final_Result}, the number of monitored lines and the lines eligible for switching change from time to time based on network loading condition. 
%\begin{frame}
%\frametitle{Case study and numerical results}
%\framesubtitle{MLL and SLL}

\begin{table}
\centering
\caption{Size of MLL and updated Switching Lines List ($SLL_u$) }
%\vspace{-0.1in}
\label{tab:Final_Result}

\begin{tabular}{| l |c| c| c| c| c|} \hline
& t=1 & t=2 & t=3 & t=4& t=5\\ \hline
	$|MLL|$ & 8 & 9 & 10 & 6 & 7\\ \hline
	$|SLL_u|$ & 18 & 11 & 3& 11 & 23\\ \hline
	
	\end{tabular}

\end{table}

%\begin{frame}
%\frametitle{Case study and numerical results}
%\framesubtitle{Impact of reducing monitoring lines and reducing switching lines liston switching plan}
%
%{\footnotesize
%%    Extra operation costs: 0.021\%\\ 
%%     Switching reduction: 244\%
% \begin{figure}
%%          \centering
%          \subfigure[]{
%          \includegraphics[width=0.37\linewidth]{SNM.png}
%      %    \caption{13 bus network with existing lines, generators and loads.}
%        }
%         \subfigure [ 
%         ]{
%            \includegraphics[width=0.37\linewidth]{SP.png}
%        %    \caption{13 bus network with existing lines, generators and loads.}
%          }
%          \caption{Switching plan without (a) and with (b) MLL and SLL }
%          \end{figure}
% 
%}
%
%\end{frame}

%\begin{frame}
%\frametitle{Case study and numerical results}
%\framesubtitle{Impact of MLL and SLL on operation costs and simulation time}
TS optimization problem is solved with $SLL_u$ instead of $SLL_o$, and the result is shown in~\figurename{~\ref{fig:SP}} (b). By comparing figures (a) and (b), it is clear that TS after applying MLL and SLL is no longer optimal, and the extra operation cost is 0.0118\% for this case. However, the simulation time is reduced from 271.67 seconds  to 2.1 seconds that shows more than 129 times simulation time reduction. Moreover, the number fo switching is reduced by 80\%. Therefore, using the proposed heuristic method is a trade-off between optimality and simulation time.

%$N-1$ contingency analysis is not integrated into the switching formulation, but as system should be operate in a way 
%\begin{itemize}
%\item Switching reduction: 80\%
%\item Extra operation costs: 0.012\%.
%\item Simulation time: reduced from 271.67 seconds to 2.1 seconds (129 times reduction).
%\item Selected switching plan with MLL will provide better property for $N-1$ contingency analysis satisfaction, as it will not let switching for lines that will cause overloads on already heavily loaded lines.
%\item For this case study, the selected switching plan with MLL and SLL  satisfies $N-1$ criterion, but the selected plan without them will not.
%\end{itemize}
% 

%\end{frame}

\subsection{Reduced ERCOT System}
A reduced model of the ERCOT system is provided in~\cite{Park_2013}. This network contains 317 buses, 427 branches, 489 conventional power plants, 36 wind farms and 182 load centers. The purpose of developing this case was to evaluate the impact of large penetration of wind in Competitive Renewable Energy Zone (CREZ) area on ERCOT market and transmission expansion requirements to transfer wind power to central and east Texas. For this reason, west Texas is simulated in detail, and the rest of the ERCOT area is aggregated to three zones as delivery points of CREZ. Parameters are set as follows:
%\begin{frame}
%\frametitle{Case study and numerical results}
%\framesubtitle{Reduced ERCOT System~\cite{Park_2013}}

\begin{itemize}
%\item Number of Buses: 317.
%\item Number of Branches: 427.
\item T=3.
\item $H1_l$=2.
\item $H2^t$=8.
\item $\alpha=50\%$.
\end{itemize}

%\end{frame}

%\begin{frame}
%\frametitle{Case study and numerical results}
%\framesubtitle{Summery of results}
MLL and SLL for this case study are shown in Table~\ref{tab:Final_Result2}. As shown in the second row, updated switching lines list $SLL_u$ contained much less lines compared to original switching lines list ($|SLL_o|=427$) that will reduce computational time significantly. 

\begin{table}
\centering
\caption{Size of MLL and updated Switching Lines List ($SLL_u$) }
%\vspace{-0.1in}
\label{tab:Final_Result2}

\begin{tabular}{| l |c| c| c| } \hline
& t=1 & t=2 & t=3 \\ \hline
	$|MLL|$ & 42 & 42 & 43\\ \hline
	$|SLL_u|$ & 91 & 97 & 63\\ \hline
	
	\end{tabular}

\end{table}	

Here is a summary of results:\\
Operation cost saving as a result of applying TS for the case w/o extra constraints on switching: 5.1\% \\
Extra cost as a result of adding new constraints: 2.1\% (3\% saving on operation costs) \\
Number of switching w/o extra constraints: 308 (22 lines are switched 3 times each)\\
Number of switching w/ extra constraints: 23 (1239\% switching reduction)\\
Simulation time w/o MLL and SLL: no answer after 2 days\\
Simulation time w/ MLL and SLL: 33 mins \\

This case study shows the impact of the proposed heuristic method on reducing the computational time. 

%\end{frame}

\section{Conclusion}\label{sec:conclude}
%\subsection{Conclusion and future work}
%\begin{frame}
%\frametitle{Conclusion and future work}
In summary:
\begin{itemize}
\item Adding some extra constraints may not significantly decrease the benefits of transmission switching, but can significantly decrease the total number of switching that benefits system reliability.
\item Considering multiple hours for TS planning with switching costs  may prevent frequent switching of a small group of lines.
\end{itemize}
As a part of our future work:
\begin{itemize}
\item Developing algorithms to decrease computational time for large scale systems 
\item Integrating network losses and reactive power requirements
\item Integrating $N-1$ contingency analysis into transmission switching as power system should be operated in a way that satisfies $N-1$ criterion. 
\end{itemize}

\section*{Acknowledgment}
The authors were supported, in part, by the Defense Threat Reduction Agency and the National Science Foundation.\\

%\end{frame}

%    \subsection{References}
%    \begin{frame}[allowframebreaks]\frametitle{References} 
%    \bibliographystyle{IEEEtran}
%    \vspace{-0.2in}
%    \scriptsize
%    \renewcommand{\baselinestretch}{0.95}
%    \bibliography{TransExpan}
%    \end{frame}

%\begin{frame}
%\centering
%Thanks for your time!\\
%Any Question?\\
%
%
%
%For further questions, please contact me: \texttt{m.majidi@utexas.edu}\\
%Slides can be downloaded from:\\
%\url{http://arxiv.org/abs/1507.03825}
%
%\end{frame}

\bibliographystyle{IEEEtran}
    \vspace{-0.2in}
    \scriptsize
    \renewcommand{\baselinestretch}{0.95}
    \bibliography{TransExpan}

\end{document}